\documentclass[letterpaper, 10pt, conference]{ieeeconf}      

\IEEEoverridecommandlockouts                              
\usepackage{graphics} 
\usepackage{epsfig} 
\usepackage{mathptmx} 
\usepackage{times} 
\usepackage{amsmath} 
\usepackage{amssymb}  

\usepackage[mathcal]{euscript} 
\usepackage{mathtools}
\usepackage{mathrsfs}

\usepackage{graphicx}
\graphicspath{{Figures/}} 
\usepackage{caption} 
\usepackage{subfig} 

\usepackage{url}

\usepackage{xcolor}
\definecolor{USred}{rgb}{0.74,0.1,0.1}
\definecolor{USblue}{rgb}{0.2,0.2,0.7}
\definecolor{green1}{cmyk}{0.82,0,1,0.1}

\definecolor{Royalblue}{cmyk}{1,0.30,0.2,0.2}

\newcommand{\numberset}{\mathbb}

\newcommand{\RR}{\numberset{R}}

\newcommand{\D}{\mathcal{D}}

\newcommand{\G}{\mathcal{G}}
\newcommand{\I}{\mathcal{I}}
\newcommand{\J}{\mathcal{J}}

\newcommand{\N}{\mathcal{N}}
\newcommand{\Q}{\mathcal{Q}}

\newcommand{\argmin}{\operatornamewithlimits{argmin}}
\newcommand{\argmax}{\operatornamewithlimits{argmax}}
\renewcommand{\vec}{\boldsymbol}
\DeclareMathOperator{\tr}{tr}

\DeclareRobustCommand{\vect}[1]{
	\ifcat#1\relax
	\boldsymbol{#1}
	\else
	\mathbf{#1}
	\fi}

\DeclareRobustCommand{\rv}[1]{\vec{\mathsf{#1}}}

\newcommand\indep{\protect\mathpalette{\protect\independenT}{\perp}}
\def\independenT#1#2{\mathrel{\rlap{$#1#2$}\mkern2mu{#1#2}}}

\newtheorem{problem}{Problem}
\newtheorem{remark}{Remark}

\title{\LARGE \bf Link Prediction: A Graphical Model Approach}

\author{Daniele Alpago, Mattia Zorzi, Augusto Ferrante 
\thanks{}
\thanks{D. Alpago, M. Zorzi and A. Ferrante are with the Department of Information Engineering, University of Padova, Padova, Italy; email:	 
	    {\tt\small alpagodani@dei.unipd.it} (D. Alpago)
        {\tt\small zorzimat@dei.unipd.it} (M. Zorzi)
        {\tt\small augusto@dei.unipd.it} (A. Ferrante)}%
\thanks{}%
}

\begin{document}

\maketitle
\thispagestyle{empty}
\pagestyle{empty}

\begin{abstract}
We consider the problem of link prediction in networks whose edge structure may vary (sufficiently slowly) over time.
This problem, with applications in many important areas including social networks, has two main variants: the first, known as {\em positive  link prediction} or PLP consists in  estimating  the appearance of a link in the network. The second, known as {\em negative link prediction} or NLP consists in  estimating  the disappearance of a link in the network.
We propose a data-driven approach to estimate the appearance/disappearance of edges. Our solution is based on a regularized optimization problem for which we prove existence and uniqueness of the optimal solution.  
    \end{abstract}


\section{Introduction}
The increased popularity of research in social-network analysis can be mainly attributed to explosion of social networks such as Facebook, Twitter, YouTube, etc. Despite the fact that one can collect plenty of information coming from social platforms, the topological and non-topological changes over time remain unknown due to the dynamical behavior exhibited by this kind of networks. Such underlying dynamical nature strongly motivates the interest in problems such as inferring unobserved edges in the current network or predicting edges that appear/disappear in the future, which are both referred to as \emph{Link Prediction} problems \cite{Wang2015,liben2007}. Link prediction techniques find important applications in e-mails networks \cite{Wang2007} and in gene-expression networks in biology \cite{NLP_Almansoori2012} just to mention a few. Generally speaking, link prediction problems can be basically distinguished in two categories: \emph{Positive Link Prediction (PLP)} and \emph{Negative Link Prediction (NLP)} problems. The most part of the research effort has been focused on the positive link prediction problem which aims to predict the formation of edges in the network. Although links' disappearance models important behaviors in social networks such as the ``\emph{unfollow}" behavior in on-line social networks, less attention has been dedicated to the NLP problem, understood as the prediction of disappearing links in the network \cite{NLP_Almansoori2012}.

Many algorithms solving the PLP problem have been proposed in the literature. In general terms, most of these methods deal with PLP over static networks and concern the construction of the so-called \emph{score matrix}, whose entry $(i,j)$ measures the probability of appearance of an edge between node $i$ and node $j$. Each element $(i,j)$ of the matrix is therefore determined by a proper choice of a \emph{similarity measure} which, given the available prior information, estimates how likely the appearance (in the near future) of an edge between two specific nodes is. These similarity measures are based on some knowledge of the topology of the network in the past and some (usually qualitative) properties the network is expected to have. For instance, a real network is expected posses the so-called ``small world'' property meaning that the network is highly clusterized and most pairs of nodes are related through short chains. As a consequence, a couple of unconnected nodes whose path length is short should have a high similarity score. Among the most successful similarity measures that have been proposed, we recall common neighbors, Adamic/Adar, Katz and spreading activation measures, see \cite{liben2007} for more details. Few similarity measures have been proposed instead in the NLP framework. These are based on a procedure that essentially ``reverses'' the PLP process based on a certain similarity measure. The method proposed in \cite{NLP_Almansoori2012}, for instance, classifies as ``likely to disappear'' those edges that (after being removed) have a small similarity measure according to the selected PLP criterium; in other words, the edges that are estimated to disappear are the ones between the nodes that have the smallest similarity measure according to a PLP algorithm applied to the network deprived of those edges.

However, the aforementioned approaches have the following limitations:
\begin{itemize}
\item there are instances called ``unfriendly prediction networks" for which most of the measures provide poor prediction performances, see for instance \cite{gao2015link}; this evidence is due by the fact that the expected network properties do not  coincide with the actual ones;
\item in many applications both   a {\em prior} for the network topology and a {\em prior} for the mathematical model are available;
moreover, the available data allow to estimate the entire model instead of only its support and the latter contains far less information than the former; however, the similarity measures in the literature only exploit the topology of the graph.
\end{itemize}

The aim of this paper is to address the above limitations. More precisely, we consider the link prediction problem over static networks modeled by means of undirected graphical models for Gaussian random vectors. Then, we propose a data-driven link prediction approach based on a similarity measure accounting the knowledge on the past information (including the topology and the mathematical model) and some noisy piece of information. The latter is represented by noisy measurements of the network in which the new links have already appeared/disappeared. Therefore: 1) we let emerge the network properties from the measurements rather than expected properties; 2) we use all the available information. It is worth remarking that it is standard in the literature to refer to a link prediction problem either if the focus is on a \emph{prediction} problem or if the problem at hand is actually a \emph{detection} problem in which the prediction comes from a combination of the past information (the current network) and some information (data) coming from the network where the new link has already appeared/disappeared. Our approach falls in the latter case. Moreover, it naturally settles into a covariance estimation problem with prior set-up \cite{5953479}: we are looking for the graphical model that agrees with the data while being as close as possible to the current graphical model (i.e. the prior). The latter can also be understood as a static version of the covariance extension problems with prior for which a large body of literature is available \cite{GeoLind2003,ringh2016,ENQVIST1, RING_2018,enqvist2001homotopy,zorzi2014,P-F-SIAM-REV,Baggio-TAC-18,Zhu-GB-TAC-19}. Our problem is formalized as an optimization problem whose solution will be proved to exist and to be unique.

Our main modeling tool will be graphical models for Gaussian random vectors which encode conditional dependence relations among a set of jointly normally distributed random variables \cite{Lauritzen1996}. The reason why we consider conditional dependence among pairs of nodes is because it  accounts for the information distributed on the whole network and has therefore a global nature which seems suitable for many of the applications mentioned before.

The topology of a Gaussian graphical model is reflected by the support of the inverse covariance matrix. Accordingly, the topology can be inferred from data by resorting to a \emph{regularized maximum likelihood} problem where the regularizer induces sparsity on the inverse covariance matrix, see \cite{friedman2008,banerjee2006}. Our link prediction approach has a similar spirit. The main differences are that: we consider the prior; we propose three different  regularizers corresponding to PLP, NLP and a mixed version of PLP and NLP, respectively. Moreover, our paradigm is strictly related to a generalized version of the Dempster's problem \cite{dempster1972}, which has been extensively studied and generalized both in the static  
\cite{6315639} and dynamic   \cite{avventi2013,ZorzSep2016,AlpZorzFer2018,zorzi2019graphical,TAC19,ZORZI2019108516,Alpago_TAC} case.

\emph{Outline of the paper.} In Section \ref{sec:GM} we review the basic properties of Gaussian graphical models that we need throughout the paper. In Section \ref{sec:LP} we formulate the link prediction problems by means of Gaussian graphical models. In Section \ref{sec:ER} some numerical experiments are presented. Finally in Section \ref{sec:CO} we draw the conclusions.

\subsection{Notation}
Given a $p\times p$ real matrix $A\in\RR^{p\times p}$, we will denote with $A^\top$ its transpose, by $A^{-1}$ its inverse, and by $\det A$ and $\tr A$ its determinant and trace, respectively. The writings $A\ge 0$ and $A>0$ denote the fact that $A$ is positive semi-definite and positive definite. The operator $\text{diag}:\RR^{p\times p}\to\RR^{p\times p}$ maps the matrix $A$ to the diagonal matrix $\text{diag}(A)$ obtained by setting to zero all the elements of $A$ outside the main diagonal. Given support a $\Omega\subseteq\{(i,j):\,i,j=1,\dots,p\}$, we denote by $\Omega^c$ the complement of $\Omega$; the map $\mathsf{P}_\Omega:\RR^{p\times p}\to\RR^{p\times p}$, defined by
\[
    [\mathsf{P}_\Omega(A)]_{ij}=
    \left\{
    \begin{split}
       &A_{ij}&\quad\text{if}\quad &(i,j)\in\Omega,\\
       &0&\quad\text{if}\quad &\text{otherwise},
    \end{split}
    \right.
\]
is the orthogonal projection onto the subspace of matrices with support $\Omega$. The set of all symmetric $p\times p$ matrices will be denoted by $\mathcal{S}_p$ while $\mathcal{S}_p^+\subset\mathcal{S}_p$ will denote the cone of the positive definite matrices in $\mathcal{S}_p$. $I_p$ denotes the identity matrix of order $p$. The symbol $\propto$ means ``proportional to''.

\section{Gaussian Graphical Models}\label{sec:GM}
Let $\rv{x}$ be an $m$-dimensional Gaussian random vector and denote with $\vec{x}_1,\dots,\vec{x}_m$ its components. Let $\G=(V,E)$ be an undirected graph with vertexes $V=\{1,\dots,m\}$ and edges $E\subseteq V\times V$. We say that the random vector $\rv{x}$ satisfies the (undirected) Gaussian graphical model $\G$, if $\rv{x}$ admits probability density function $\N(0,\Sigma)$ with $\Sigma>0$ and such that
\begin{equation}\label{eq:supp}
	(\Sigma^{-1})_{ij}=0,\qquad\text{for all}\qquad (i,j)\notin E.
\end{equation}
Accordingly, a complete characterization of the graph $\G$ associated to the random vector $\rv{x}$ is given in terms of the support of its \emph{concentration matrix} $\Sigma^{-1}$. One can show \cite{Lauritzen1996} that the random variable $\vec{x}_i$ is conditional independent from $\vec{x}_j$ given the other components $\{\vec{x}_k\}_{k\ne{i,j}}$, and we write $\vec{x}_i \indep \vec{x}_j \,|\, \{\vec{x}_k\}_{k\ne i,j}$, if and only if $(\Sigma^{-1})_{ij}=0$. Therefore,
\begin{equation}\label{eq:ci}
	\vec{x}_i \indep \vec{x}_j \,|\, \{\vec{x}_k\}_{k\ne i,j}
	\quad\iff\quad
	(i,j)\notin E.
\end{equation} The identification of Gaussian graphical models from data, in view of \eqref{eq:ci}, reduces to the estimation of the covariance matrix $\Sigma$ whose inverse should be sparse. In particular, one can use the maximum likelihood principle. Given $N$ i.i.d. observations $\mathsf{x}_1,\dots,\mathsf{x}_N$ of $\rv{x}$, the negative log-likelihood function has the form
\begin{equation}\label{eq:like}
	\ell(\Sigma; \hat \Sigma) \propto \frac{N}{2}\log\det\Sigma + \frac{N}{2}\tr\left(\hat{\Sigma}\,\Sigma^{-1}\right)
\end{equation}
where
\begin{equation}\label{eq:samplecov}
	\hat{\Sigma} = \frac{1}{N}\sum_{k=1}^{N}\,\mathsf{x}_k\,\mathsf{x}_k^\top
\end{equation}
is the sample covariance matrix. Hence, the maximum likelihood estimator of $\Sigma$ boils down to the constrained optimization problem
\begin{equation}\label{eq:maxlik}
	\argmin_{\Sigma\in\Theta} \quad \log\det\Sigma + \tr\left(\hat{\Sigma}\,\Sigma^{-1}\right)
\end{equation} 
where $\Theta\subseteq\mathcal{S}_p^+$ is a suitable parametric family (e.g. the family of covariance matrices corresponding to a graphical model having a certain topology).

\section{Link Prediction}\label{sec:LP}
This section is devoted to the introduction and the formalization of the link prediction problem in our setting. Let $\G_s=(V,E_s)$ be the graphical model associated to an $m$-dimensional Gaussian random vector $\rv{x}$ with covariance matrix $S$, modeling some kind of network at a certain time $s>0$. By relation \eqref{eq:supp} the support of $S^{-1}$,
\[
\Omega_s:=\left\{(i,j)\in V\times V:\,(S^{-1})_{ij}\ne0\right\},
\]
coincides with the set $E_s$ and it is assumed to be \emph{known}. 
Adopting the standard assumption that the number of nodes does not change over time, we model our system at time $t>s$ with the graphical model $\G_t=(V,E_t)$ associated to the same random vector $\rv{x}$. The (new) edges' set is now related to the concentration matrix $T^{-1}$, through its support $\Omega_t$ which is considered \emph{uknown}, and therefore it has to be estimated.

Given the application, it is reasonable to assume that $T$ will not be drastically different from $S$ because this would mean that the network has completely changed in a relatively short period of time, which is typically not the case. 
Notice that this is not a formal mathematical assumption but must be regarded simply as a justification of the fact that among all the models compatible with the data we select the one that is closer to to the `prior' given by $S$.
We are now ready to state the following estimation problem.

\begin{problem}\label{pb:linkcompl}
Assume to have the ``prior'' covariance matrix $S$ of $\rv x$ at time $s$. Let $\Omega_s$ be the support of $S^{-1}$. Given $N$ i.i.d. observations $\mathsf{x}_t^1,\dots,\mathsf{x}_t^N$ of $\rv{x}$ at time $t$, with $t>s$, compute an estimate of the covariance $T$, which is close to $S$ and fits the observation as much as possible.\\ 
\end{problem}

Notice that, Problem \ref{pb:linkcompl} includes the main variations of the link prediction problem. In fact, in link prediction we only need to estimate the support $\Omega_t$ (which is fully specified by $T$).
In particular, 
\begin{itemize}
	\item[-] $\Omega_t \supset \Omega_s$ corresponds to the PLP problem since conditional dependencies, and hence edges, are appearing between the variables;
	\item[-] $\Omega_t \subset \Omega_s$ corresponds to the NLP problem since edges are disappearing between the variables.	
\end{itemize}

\begin{figure}[h!]\centering
	\includegraphics[scale=0.8]{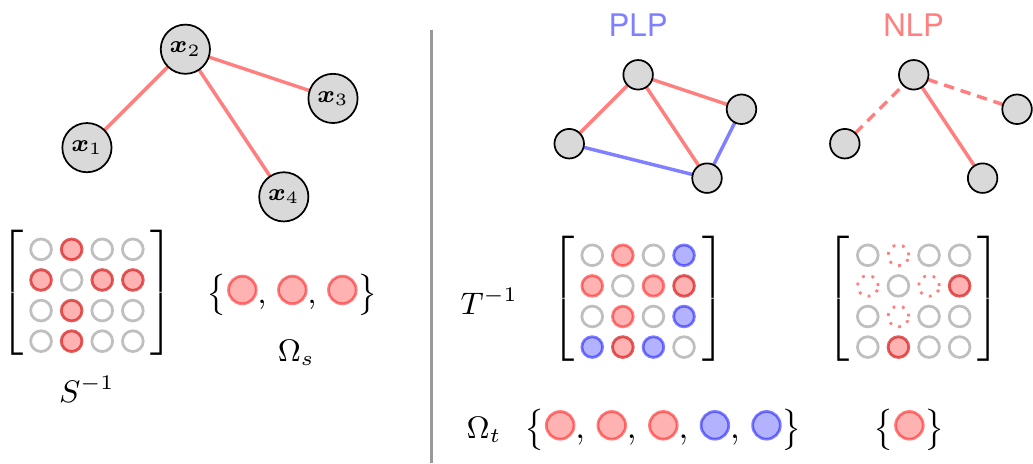}
	\caption{Link prediction problems for a four-nodes graph. The dashed red edges disappears at time $t$ while the blue ones appears. For simplicity $\Omega_s$ and $\Omega_t$ are identified with half of the support of $S^{-1}$ and $T^{-1}$, respectively.}\label{fig:setup}
\end{figure}

Figure \ref{fig:setup} explains the problems just presented for a simple network of four nodes.

We first make the simplistic assumption that we know the support $\Omega$. The solution to Problem \ref{pb:linkcompl} that we will derive under this assumption, will give us insights on how to deal with the actual problem where $\Omega$ is unknown. Let $\hat{T}$ be the sample covariance matrix \eqref{eq:samplecov} computed with the observations $\mathsf{x}_t^1,\dots,\mathsf{x}_t^N$.  In light of the previous observations, we propose the following mathematical formulation of Problem \ref{pb:linkcompl}:
\begin{equation}\label{eq:KLmin}
	\begin{aligned}
		\argmin_{T\in\mathcal{S}_p^+} &\quad 2 \D(T\|S)\\
		\text{subject to } &\quad \mathsf{P}_\Omega\left(T - \hat{T}\right)=0,
	\end{aligned}
\end{equation}
where
\begin{equation}\label{eq:KL}
	\D(T\|S):=\frac{1}{2}[-\log\det (S^{-1}T)+\tr\left(S^{-1}T\right)-m]
\end{equation} 
is the Kullback-Leibler divergence between the distributions $\N(0,S)$ and $\N(0,T)$. Some remarks concerning Problem \eqref{eq:KLmin} are in order. First of all, notice that if $S=I$ and $\{(i,i):\,i=1,\dots,m\}\subset\Omega_t$, Problem \eqref{eq:KLmin} becomes the classical Dempster's problem, widely studied in the literature \cite{dempster1972}. Indeed, in this case, the linear term $\tr(S^{-1}T)=\tr T$ is in fact the constant $\tr\hat{T}$, thanks to the constraint.
The index that has to be optimized in \eqref{eq:KLmin} imposes invertibility of the solution $T$ so that even if the ``true'' covariance is singular we find a nonsingular approximation.

The constraint in \eqref{eq:KLmin} imposes that the entries of $T$ coincide with those of  $\hat{T}$ on the support $\Omega$. Enforcing such a constraint may be confusing at first sight, as the edges of the predicted network have to do with the support of the \emph{concentration} matrix rather than the support of the \emph{covariance} matrix. However, conditional dependence relations among the variables are strictly related to their correlations and with that constraint we are precisely accounting \emph{only} the most reliable correlations inferred from the data. As we will see at the end of the section, such constraint imposes the right structure of the solution $T_o^{-1}$ which will have as support the union $\Omega_s\cup\Omega$, that will allow to model the different link prediction scenarios.

The problem can be re-parametrized in term of $T^{-1}$ exploiting duality theory. 
To this end, we first eliminate the uninteresting constant term $\log\det(S) -m$ form the cost function in \eqref{eq:KLmin}.
Then we form the Lagrangian for Problem \eqref{eq:KLmin}: it is
\begin{align*}
	\mathcal{L}(T,\tilde{\Lambda})&=-\log\det T+\tr\left(S^{-1} T\right)+\tr\left[\mathsf{P}_{\Omega}\left( T - \hat{ T}\right)\,\tilde{\Lambda}\right]\\
	&=-\log\det T+\tr\left(S^{-1} T\right) +\tr\left[\left( T - \hat{ T}\right)\,\mathsf{P}_{\Omega}\left(\tilde{\Lambda}\right)\right]
\end{align*} 
where $\tilde{\Lambda}=\tilde{\Lambda}^\top\in\RR^{m\times m}$ is the Lagrange multiplier and we have exploited the fact that $\mathsf{P}_{\Omega}(\cdot)$,
being an orthogonal projection, is  a self-adjoint operator as it may be readily checked by applying the definition. Introducing the new multiplier $\Lambda:=\mathsf{P}_{\Omega}(\tilde{\Lambda})$, we can rewrite the Lagrangian as
\begin{equation*}
	\mathcal{L}( T,\Lambda)=-\log\det T+\tr\left[(S^{-1}+\Lambda) T\right]-\tr\left(\hat{ T}\,\Lambda\right).
\end{equation*} 
Since $\mathcal{L}$ is a strictly convex function of $ T$, a sufficient condition for $ T_o$ to be a minimum point for $\mathcal{L}$ is that the first variation of $\mathcal{L}$ in direction $\delta T$ is zero for every direction $\delta T\in\mathcal{S}_p$, namely
\begin{equation*}
	\mathcal{L}( T,\Lambda;\delta T)=\tr(- T^{-1}\delta T+\left(S^{-1}+\Lambda\right)\delta T)=0,\,\forall\,\delta T\in\mathcal{S}_p.
\end{equation*}
The form of the minimum is therefore $ T_o=\left(S^{-1}+\Lambda\right)^{-1}$ provided that $\Lambda\in\mathcal{S}_p$ is chosen so that $S^{-1}+\Lambda>0$. The dual of Problem \eqref{eq:KLmin} is therefore 
\begin{equation}\label{eq:dualpb-vero}
	\begin{aligned}
		\argmax_{\Lambda\in\Q_S} &\quad \mathcal{L}( T_o,\Lambda)\\
		\text{subject to } &\quad \mathsf{P}_{\Omega^c}\left(\Lambda\right)=0
	\end{aligned}
	\end{equation}
and $\Q_S:=\left\{\Lambda\in\mathcal{S}_p:\,S^{-1}+\Lambda>0\right\}$ is the domain of optimization. To remain in the convex setting, we change sign and minimize the opposite function; namely, we consider the following problem:
\begin{equation}\label{eq:dualpb}
	\begin{aligned}
		\argmin_{\Lambda\in\Q_S} &\quad \J_S(\Lambda)\\
		\text{subject to } &\quad \mathsf{P}_{\Omega^c}\left(\Lambda\right)=0,
	\end{aligned}
\end{equation}
where the dual functional (save for constant terms) is
\begin{equation*}
	\J_S(\Lambda)=-\mathcal{L}( T_o,\Lambda) +m=-\log\det\left(S^{-1}+\Lambda\right)+\tr\left(\hat{ T}\,\Lambda\right).
\end{equation*}

\begin{remark}
Performing the substitution $T=(S^{-1}+\Lambda)^{-1}$ in (\ref{eq:dualpb}), the dual functional $\J_S$ correspond (save for constant factors) to the negative log-likelihood $\ell( T;\hat{ T})$ in \eqref{eq:like}. Hence, in view of \eqref{eq:maxlik}, Problem \eqref{eq:dualpb} can be readily interpreted as a maximum-likelihood problem,
	\begin{equation}\label{eq:MLdual}
		\argmin_{ T\in\Theta_{S,\Lambda}} \quad \ell( T;\hat{ T}),
	\end{equation}	
	over the parametric family
	\[
	\Theta_{S,\Lambda}:=\left\{  T=(S^{-1}+\Lambda)^{-1}:\,\Lambda\in\Q_S,\,\mathsf{P}_{\Omega^c}\left(\Lambda\right)=0 \right\}.
	\]
\end{remark}

\vspace{3mm}

The derivation of the dual problem shows that the support of the optimum $ T_o^{-1}$ is given by $\Omega_s\cup\Omega$. If $\Omega_t$  were known and if the relation $\Omega_t \supset \Omega_s$ holds, then we  could choose $\Omega=\Omega_t$ so that the support of $ T_o^{-1}$ would be $\Omega_t$;  in this case, we have a PLP problem. Similarly, knowing $\Omega_t$ and if the relation $\Omega_s \supset \Omega_t$ holds, then we may choose $\Omega=\Omega_s$ so that the support of $ T_o^{-1}$ would be $\Omega_t$, \emph{provided that}  $\Lambda$ is such that $(\Lambda)_{ij}=-(S^{-1})_{ij}$, $\forall\, (i,j)\in \Omega_s \setminus \Omega_t$: this is the NLP problem. The requirement on $\Lambda_{ij}$ over $\Omega_s\setminus\Omega_t$ produces zeros on $T^{-1}$ so that its support becomes smaller than $\Omega_s$, which is precisely what happens in NLP. In practice $\Omega_t$ is \emph{unknown} and the aforementioned constraints cannot be enforced. Equivalently, in view of \eqref{eq:MLdual}, the parametric family $\Theta_{S,\Lambda}$ is unknown.
To overcome this issue, we estimate $\Omega_t$ from the available data at time $t$ by considering two different regularized versions of Problem \eqref{eq:MLdual} corresponding to the two link prediction problems:
\begin{itemize}
	\item[-] The PLP problem, for which $\Omega_t\supset\Omega_s$, boils down to the regularized maximum-likelihood problem
	\begin{equation}\label{eq:dualPLP}
		\argmin_{ T\in\Theta_{S,\Lambda}^P} \quad \ell( T;\hat{ T})+\gamma_P\,h_P( T),
	\end{equation}
	where $\Theta_{S,\Lambda}^P:=\left\{  T=(S^{-1}+\Lambda)^{-1}:\,\Lambda\in\Q_S\right\}$ and 
	\begin{equation}\label{eq:l1PLP}
	h_P( T)=\sum_{(i,j)\in\I_P}\,|\Lambda_{ij}|
	\end{equation}
	with $\I_P := \left\{(i,j)\in V\times V\setminus\Omega_s:\,i>j\right\}$. With this choice we are inducing sparsity among the $\Lambda_{ij}$'s \emph{outside} $\Omega_s$ so that the support $\Omega_t$ of $ T_o^{-1}=S^{-1}+\Lambda$ will contain $\Omega_s$.
	It is worth noticing that from the definition of the index sets $\I_P$ we are not penalizing the entries fixed by the prior support $\Omega_s$. This will reduce the bias (i.e. the shrinking to zero of the entries) affecting the final estimate.
	\item[-] The NLP problem, for which $\Omega_t\subset\Omega_s$, can be formulated as a regularized maximum-likelihood problem as well,
	\begin{equation}\label{eq:dualNLP}
		\argmin_{ T\in\Theta_{S,\Lambda}^N} \quad \ell( T;\hat{ T})+\gamma_N\,h_N( T),
	\end{equation}
	where $\Theta_{S,\Lambda}^N:=\left\{  T=(S^{-1}+\Lambda)^{-1}:\,\Lambda\in\Q_S,\,\mathsf{P}_{\Omega_s^c}(\Lambda)=0\right\}$
	
	 and 
	\begin{equation}\label{eq:l1NLP}
	h_N( T)=\sum_{(i,j)\in\I_N}\,|\Lambda_{ij}+(S^{-1})_{ij}|
	\end{equation}
	with $\I_N := \left\{(i,j)\in\Omega_s:\,i>j\right\}$. Here, we force $\Lambda_{ij}=0$ on $\Omega_s^c$ and induce $\Lambda_{ij}$ to be equal to $-S^{-1}_{ij}$ \emph{over} $\Omega_s$. Hence, in this case, $\Omega_t$ will be contained in $\Omega_s$.
\end{itemize}
Here $\gamma_P,\,\gamma_N>0$ are the two regularization parameters, weighting the effect of the $\ell_1$-penalties $h_P$ and $h_N$, respectively.\\

\begin{remark}
    The approach outlined above, can be adapted in the case of a \emph{mixed link-prediction} situation, in which one intends to predict both appearing and disappearing links. This may be achieved by considering the PLP regularized maximum-likelihood problem \eqref{eq:dualPLP} and substituting the regularizer with a combination of \eqref{eq:l1PLP} and \eqref{eq:l1NLP}, 
    \begin{equation*}
    	h_M(T):=\eta_P\sum_{(i,j)\in\I_P}\,|\Lambda_{ij}|+\eta_N\,\sum_{(i,j)\in\I_N}|\Lambda_{ij}+(S^{-1})_{ij}|
    \end{equation*}
    where $\eta_P>0$ and $\eta_N>0$ are the two regularization parameters that contribute to enforce the sparsity of $\Lambda$ over $\Omega_s^c$ for the PLP and the fact that $\Lambda_{ij}=-(S^{-1})_{ij}$ over $\Omega_s$ for the NLP, respectively. The resulting regularized maximum-likelihood problem is 
    \begin{equation}\label{eq:dualMLP}
    \argmin_{ T\in\Theta_{S,\Lambda}^P} \quad \ell( T;\hat{ T})+\eta_P\,h_P(T)+\eta_N\,h_N(T).
    \end{equation}
\end{remark}

Tuning $\eta_P$ and $\eta_N$ in (\ref{eq:dualMLP}) it is possible to emphasize the PLP or the NLP task, respectively.

It can be shown that Problems \eqref{eq:dualPLP}, \eqref{eq:dualNLP} and  \eqref{eq:dualMLP} admit a unique solution.
While uniqueness follows by showing that (in an equivalent formulation of the problem) the functional is strictly convex, 
the argument for existence is based on the fact that the functional goes to $+\infty$
whenever its  argument  tends to the boundary of the 
feasible set or diverges.

\subsection{Induced Similarity Measure}
Let $\hat T_o$ be the solution to Problem (\ref{eq:dualPLP}) (respectively (\ref{eq:dualNLP}), (\ref{eq:dualMLP}), depending on the link prediction problem we are cosidering). The latter not only characterizes the graphical model $\mathcal G_t$ of $\vec x$ at time $t$, but it also induces a similarity measure between the nodes of  $\mathcal G_t$. 
 In fact, our estimation method naturally induces the \emph{score matrix} which evaluates the inclination of pairs of nodes to be conditionally dependent:
\begin{equation*}
	R_t := \text{diag}(\hat T_o)^{1/2}\, \hat T_o^{-1}\,\text{diag}(\hat T_o)^{1/2}.
\end{equation*}    
The measure of the similarity between node $i$ and node $j$ that comes out from our approach is therefore $r_{ij}=(R_t)_{ij}$. More precisely $-1\leq r_{ij} \leq 1$ and the higher $|r_{ij}|$ is  the more probable is that a link connecting node $i$ to node $j$ will appear (i.e. node $i$ and node $j$ will be conditionally dependent).It is worth noting that in order for our similarity measure to be a suitable topology-selection measure according to relations \eqref{eq:supp}-\eqref{eq:ci}, a thresholding procedure is needed:
\begin{equation}\label{eq:thresholding}
	(i,j)\in\hat{\Omega}_t \quad\iff\quad |r_{ij}|>t_r,
\end{equation}
where $t_r>0$ has to be set properly.

\section{Experimental Results}\label{sec:ER}
\begin{figure}[h!]\centering
	\includegraphics[scale=0.75]{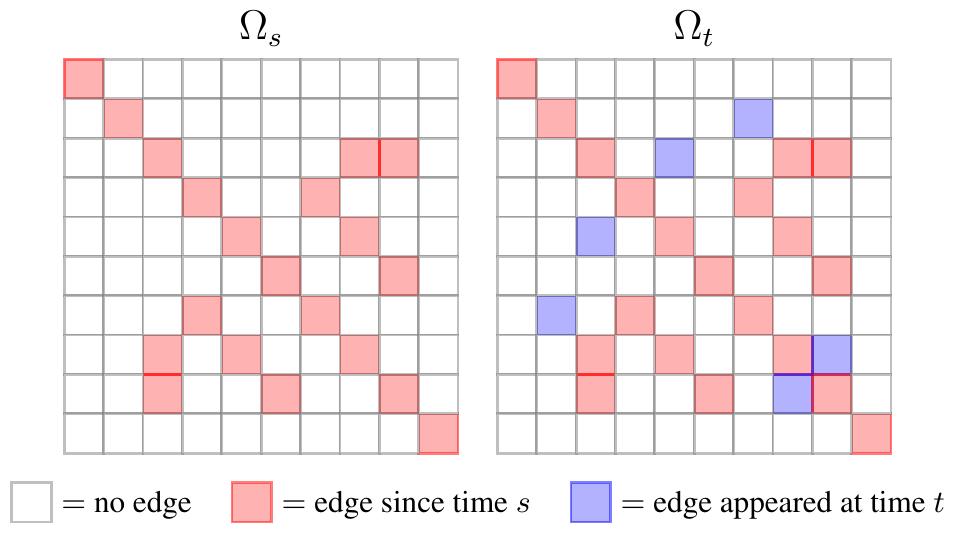}
	\caption{PLP case. Support $\Omega_s$ of the prior concentration matrix (left) and support $\Omega_t$ of the concentration matrix $T^{-1}$ that has to be estimated (right).}\label{fig:true_supp_PLP}
\end{figure}

In this section we show how the proposed method works in predicting appearing/disappearing edges in a simple network of $m=10$ nodes. The thresholding value $t_r$ is set to $10^{-4}$ throughout all the experiments. First we present the results for PLP. Figure \ref{fig:true_supp_PLP} reports the support $\Omega_s$ of the prior $S^{-1}$ at time $s$ and the true support $\Omega_t$ of the concentration matrix $T^{-1}$ at time time $t>s$, that we want to estimate.
Such a network is an unfriendly positive link prediction network. Indeed,  consider the ``common neighbors'' similarity measure \cite{liben2007} at time $t$:
\begin{align}\label{PLPSIMMEAS}
(CN_0)_{ij}=\mathrm{card}(\mathcal N_i \cap\mathcal N_j )
\end{align}
where $\mathcal N_i$ is the set of neighbors of node $i$, according to the  topology of the ``prior'' network $\mathcal G_s(V,\Omega_s)$, and  $\mathrm{card}(\mathcal N_i \cap\mathcal N_j )$ denotes the cardinality of set $\mathcal N_i \cap\mathcal N_j $. Then, it is not difficult to see that: $(CN_0)_{ij}=1$ for the pairs $(3,6)$ and $(8,9)$;  $(CN_0)_{ij}=0$ for the other $(i,j)\notin\Omega_s$. Therefore, this similarity measure  mispredicted the appearance of  3 edges.  

\begin{figure}[h!]\centering
	\includegraphics[scale=0.65]{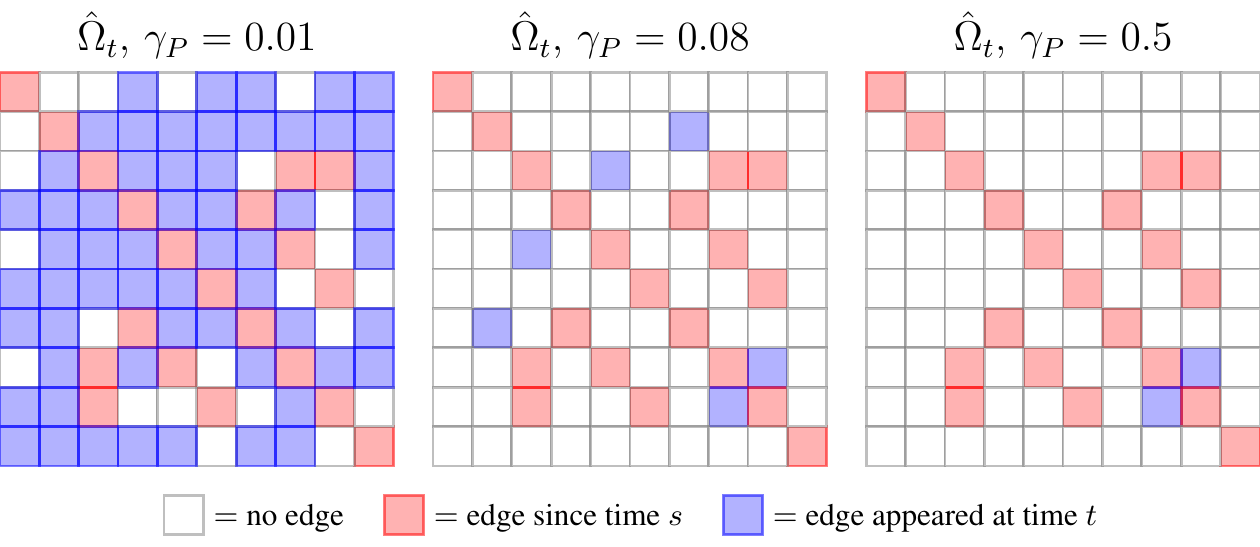}
	\caption{Estimates of the support $\Omega_t$ for different values of the regularization parameters $\gamma_P$.}\label{fig:est_supp_PLP}
\end{figure}

Starting from the left of Figure \ref{fig:est_supp_PLP} we have three versions of the estimated support $\hat{\Omega}_t$ for $\gamma_P=0.01,\,0.08,\,0.5$, respectively using our method. By comparison with Figure \ref{fig:true_supp_PLP}, $\gamma_P=0.08$ is the best choice, in that the estimation procedure infers exactly the true $\Omega_t$. For $\gamma_P=0.01$ the effect of regularization is too mild, i.e. sparsity is not properly enforced, while $\gamma_P=0.5$ gives an exaggerated sparse estimate, as one should expect. The estimation has been performed by solving Problem \eqref{eq:dualPLP} with the CVX package for Matlab \cite{GrantBoyd08,GrantBoyd14}. In particular, the sampling covariance matrix $\hat{\Sigma}$ has been computed as in \eqref{eq:samplecov} using $N=1000$ i.i.d. observations drawn from the distribution $\N(0,T)$.
\vspace{1mm}%

\begin{figure}[h!]\centering
	\includegraphics[scale=0.75]{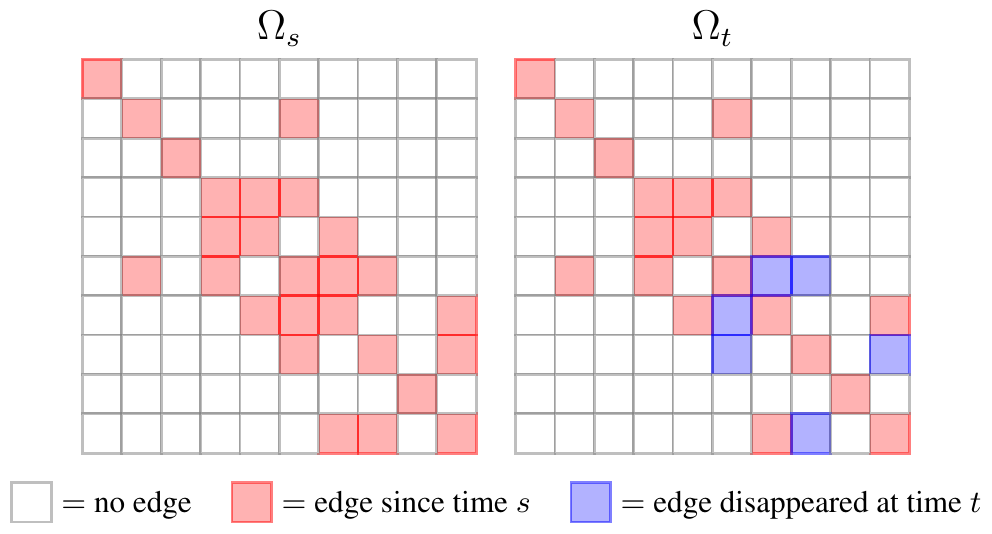}
	\caption{NLP case. Support $\Omega_s$ of the prior concentration matrix (left) and support $\Omega_t$ of the concentration matrix $T^{-1}$ that has to be estimated (right).}\label{fig:true_supp_NLP}
\end{figure}

\begin{figure}[h!]\centering
	\includegraphics[scale=0.62]{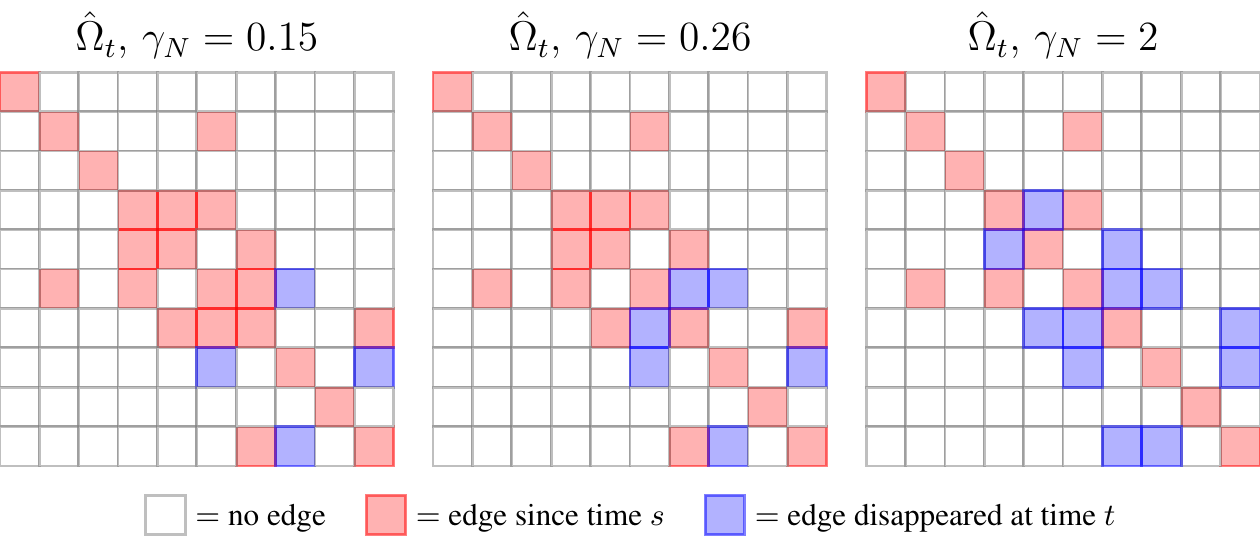}
	\caption{Estimates of the support $\Omega_t$ for different values of the regularization parameters $\gamma_N$.}\label{fig:est_supp_NLP}
\end{figure}

The NLP experiment has been performed along the same lines of the PLP one.   As before, we consider a network of $m=10$ nodes whose edges are prescribed by the support $\Omega_s$ of the prior concentration matrix $S^{-1}$, depicted in Figure \ref{fig:true_supp_NLP} (left). The true support $\Omega_t$ of the concentration matrix $T^{-1}$ has been obtained by setting to zero some of the elements in $\Omega_s$ (in blue) and it is represented in Figure \ref{fig:true_supp_NLP} (right). This network appears to have the ``unfriendly prediction'' property as the method in \cite[Sect. 4.2.2]{NLP_Almansoori2012} applied by ``reversely'' implementing the  PLP similarity measure (\ref{PLPSIMMEAS}) provides poor performances. More precisely, it correctly predicts only the disappearance of the edge (8,10), while it erroneously  predicts the disappearance of edge (4,5) and does not predict the disappearance of the two edges 
(6,7) and (6,8) that indeed disappear. Thus, in total, 3 edges are mispredicted.
\\%
Using our method, instead, the estimate $\hat{\Omega}_t$ of $\Omega_t$ has been obtained by solving Problem \eqref{eq:dualNLP} with the CVX Matlab's package. Figure \ref{fig:est_supp_NLP} depicts three estimates of $\hat{\Omega}_t$ corresponding to three values of the regularization parameters $\gamma_N=0.15,\,0.26,\,2$.
\\%
Comparing Figure \ref{fig:est_supp_NLP} and Figure \ref{fig:true_supp_NLP}, one can see that $\gamma_N=0.26$ is indeed a good choice of the regularization parameter, since produces the estimate $\hat{\Omega}_t=\Omega_t$. 

\begin{figure}[h!]\centering
	\includegraphics[scale=0.8]{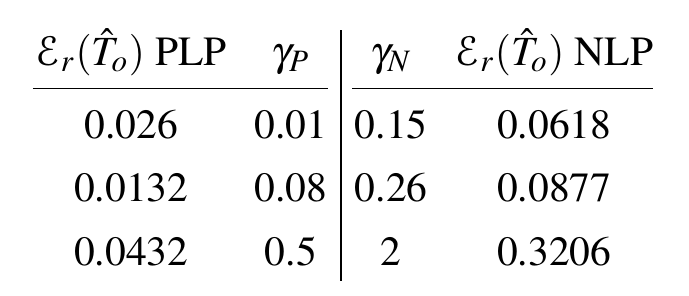}
	\caption{Relative estimation error $\mathcal{E}_r$ for the different regularization parameters $\gamma_P$s and $\gamma_S$s considered in the above PLP and NLP examples.}\label{fig:est_err}
\end{figure}

Finally, in Figure \ref{fig:est_err} we have reported the estimation errors 
$\mathcal{E}_r(\hat{T}_o) = \|T-\hat{T}_o\|_F /\|T\|_F$
for the different values of the regularization parameters $\gamma_P$s and $\gamma_N$s considered in the previous examples. For the PLP case, the model with smallest  $\mathcal{E}_r$ is the one with $\gamma_P=0.08$, as expected. In regard to the NLP case, the model with smallest  $\mathcal{E}_r$ is the one with $\gamma_N=0.15$. It is worth noting that the latter is still better than the one obtained by the procedure in \cite[Sect. 4.2.2]{NLP_Almansoori2012}; indeed, our model with $\gamma_N=0.15$ mispredicts only 1 edge.
\\%

\section{Conclusions}\label{sec:CO}
In this work an estimation method based on $\ell_1$-regularized maximum likelihood has been proposed and applied to link prediction problems. Both positive and negative link prediction problems have been formulated as convex optimization problems, whose solution is unique. The most significant contribution of the work is the introduction of a similarity measure exploiting not only the topology of the current network (prior), but also noisy information coming from the network at the current time rather than some properties the network is expected to have. The numerical examples have shown that our method can be a potential alternative to some of the existent link prediction approaches.


\begin{thebibliography}{10}

\bibitem{Wang2015}
P.~Wang, B.~Xu, Y.~Wu, and X.~Zhou, ``Link prediction in social networks: the
  state-of-the-art,'' {\em Science China Information Sciences}, vol.~58,
  pp.~1--38, Jan 2015.

\bibitem{liben2007}
D.~Liben-Nowell and J.~Kleinberg, ``The link-prediction problem for social
  networks,'' {\em Journal of the American society for information science and
  technology}, vol.~58, no.~7, pp.~1019--1031, 2007.

\bibitem{Wang2007}
C.~Wang, V.~Satuluri, and S.~Parthasarathy, ``Local probabilistic models for
  link prediction,'' in {\em Proceedings of the 2007 Seventh IEEE International
  Conference on Data Mining}, ICDM '07, (Washington, DC, USA), pp.~322--331,
  IEEE Computer Society, 2007.

\bibitem{NLP_Almansoori2012}
W.~Almansoori, S.~Gao, T.~Jarada, A.~Elsheikh, A.~Murshed, J.~Jida, R.~Alhajj,
  and J.~Rokne, ``Link prediction and classification in social networks and its
  application in healthcare and systems biology,'' {\em Network Mod. An. Health
  Inf. and Bio.}, vol.~1, pp.~27--36, Jun 2012.

\bibitem{gao2015link}
F.~Gao, K.~Musial, C.~Cooper, and S.~Tsoka, ``Link prediction methods and their
  accuracy for different social networks and network metrics,'' {\em Scientific
  programming}, vol.~2015, p.~1, 2015.

\bibitem{5953479}
A.~{Ferrante}, M.~{Pavon}, and M.~{Zorzi}, ``A maximum entropy enhancement for
  a family of high-resolution spectral estimators,'' {\em IEEE Transactions on
  Automatic Control}, vol.~57, pp.~318--329, Feb 2012.

\bibitem{GeoLind2003}
T.~T. Georgiou and A.~Lindquist, ``Kullback-leibler approximation of spectral
  density functions,'' {\em IEEE Transactions on Information Theory}, vol.~49,
  no.~11, pp.~2910--2917, 2003.

\bibitem{ringh2016}
A.~Ringh, J.~Karlsson, and A.~Lindquist, ``Multidimensional rational covariance
  extension with applications to spectral estimation and image compression,''
  {\em SIAM Journal on Control and Optimization}, vol.~54, no.~4,
  pp.~1950--1982, 2016.

\bibitem{ENQVIST1}
C.~I. {Ryrnes}, P.~{Enqvist}, and A.~{Lindquist}, ``Cepstral coefficients,
  covariance lags, and pole-zero models for finite data strings,'' {\em IEEE
  Transactions on Signal Processing}, vol.~49, pp.~677--693, April 2001.

\bibitem{RING_2018}
A.~Ringh, J.~Karlsson, and A.~Lindquist, ``Multidimensional rational covariance
  extension with approximate covariance matching,'' {\em SIAM Journal on
  Control and Optimization}, vol.~56, no.~2, pp.~913--944, 2018.

\bibitem{enqvist2001homotopy}
P.~Enqvist, ``A homotopy approach to rational covariance extension with degree
  constraint,'' {\em International Journal of Applied Mathematics and Computer
  Science}, vol.~11, pp.~1173--1201, 2001.

\bibitem{zorzi2014}
M.~Zorzi, ``A new family of high-resolution multivariate spectral estimators,''
  {\em IEEE Transactions on Automatic Control}, vol.~59, no.~4, pp.~892--904,
  2014.

\bibitem{P-F-SIAM-REV}
M.~Pavon and A.~Ferrante, ``On the geometry of maximum entropy problems,'' {\em
  SIAM Review}, vol.~55, no.~3, pp.~415--439, 2013.

\bibitem{Baggio-TAC-18}
G.~Baggio, ``Further results on the convergence of the {P}avon-{F}errante
  algorithm for spectral estimation,'' {\em IEEE Trans. Autom. Control},
  vol.~63, pp.~3510--3515, Oct 2018.

\bibitem{Zhu-GB-TAC-19}
B.~Zhu and G.~Baggio, ``On the existence of a solution to a spectral estimation
  problem {\`a} la {B}yrnes-{G}eorgiou-{L}indquist,'' {\em IEEE Trans. on
  Autom. Control}, vol.~64, pp.~820--825, Feb 2019.

\bibitem{Lauritzen1996}
S.~L. Lauritzen, {\em Graphical models}, vol.~17.
\newblock Clarendon Press, 1996.

\bibitem{friedman2008}
J.~Friedman, T.~Hastie, and R.~Tibshirani, ``Sparse inverse covariance
  estimation with the graphical lasso,'' {\em Biostatistics}, vol.~9, no.~3,
  pp.~432--441, 2008.

\bibitem{banerjee2006}
O.~Banerjee, L.~E. Ghaoui, A.~d'Aspremont, and G.~Natsoulis, ``Convex
  optimization techniques for fitting sparse gaussian graphical models,'' in
  {\em Proceedings of the 23rd international conference on Machine learning},
  pp.~89--96, ACM, 2006.

\bibitem{dempster1972}
A.~P. Dempster, ``Covariance selection,'' {\em Biometrics}, pp.~157--175, 1972.

\bibitem{6315639}
L.~{Ning}, X.~{Jiang}, and T.~{Georgiou}, ``Geometric methods for structured
  covariance estimation,'' in {\em ACC}, pp.~1877--1882, 2012.

\bibitem{avventi2013}
E.~Avventi, A.~G. Lindquist, and B.~Wahlberg, ``{ARMA} identification of
  graphical models,'' {\em IEEE Transactions on Automatic Control}, vol.~58,
  no.~5, pp.~1167--1178, 2013.

\bibitem{ZorzSep2016}
M.~Zorzi and R.~Sepulchre, ``{AR} identification of latent-variable graphical
  models,'' {\em IEEE Transactions on Automatic Control}, vol.~61,
  pp.~2327--2340, Sept 2016.

\bibitem{AlpZorzFer2018}
D.~Alpago, M.~Zorzi, and A.~Ferrante, ``Identification of sparse reciprocal
  graphical models,'' {\em IEEE Control Systems Letters}, vol.~2, pp.~659--664,
  Oct 2018.

\bibitem{zorzi2019graphical}
M.~Zorzi, ``Graphical model selection for a particular class of continuous-time
  processes,'' {\em Kybernetika}, vol.~55, no.~5, pp.~782--801, 2019.

\bibitem{TAC19}
V.~Ciccone, A.~Ferrante, and M.~Zorzi, ``Learning latent variable dynamic
  graphical models by confidence sets selection,'' {\em IEEE Trans. Autom.
  Control (accepted)}, 2020.

\bibitem{ZORZI2019108516}
M.~Zorzi, ``Empirical {B}ayesian learning in {AR} graphical models,'' {\em
  Automatica}, vol.~109, p.~108516, 2019.

\bibitem{Alpago_TAC}
D.~Alpago, M.~Zorzi, and A.~Ferrante, ``A scalable strategy for the
  identification of latent-variable graphical models,'' {\em Submitted}, 2018.

\bibitem{GrantBoyd08}
M.~Grant and S.~Boyd, ``Graph implementations for nonsmooth convex programs,''
  in {\em Recent Advances in Learning and Control} (V.~Blondel, S.~Boyd, and
  H.~Kimura, eds.), Lecture Notes in Control and Information Sciences,
  pp.~95--110, Springer-Verlag Limited, 2008.
\newblock \url{http://stanford.edu/~boyd/graph_dcp.html}.

\bibitem{GrantBoyd14}
M.~Grant and S.~Boyd, ``{CVX}: Matlab software for disciplined convex
  programming, version 2.1.'' \url{http://cvxr.com/cvx}, Mar. 2014.

\end{thebibliography}
\end{document}